\pdfoutput=1
\RequirePackage{ifpdf}
\ifpdf 
\documentclass[pdftex]{sigma}
\else
\documentclass{sigma}
\fi

\usepackage{mathrsfs}

\def\cA{\mathcal A}
\def\cC{\mathcal C}
\def\cH{\mathcal H}
\def\cD{\mathcal D}
\def\cP{\mathcal P}
\def\diag{\hbox{diag}}
\def\hc{\hat c}

\numberwithin{equation}{section}
\newtheorem{Theorem}{Theorem}[section]
\newtheorem*{Theorem*}{Theorem}

 { \theoremstyle{definition}

\newtheorem{Example}[Theorem]{Example}
 }

\begin{document}
\allowdisplaybreaks

\newcommand{\arXivNumber}{2410.14068}

\renewcommand{\PaperNumber}{083}

\FirstPageHeading

\ShortArticleName{$q$-Hypergeometric Orthogonal Polynomials with $q=-1$}

\ArticleName{$\boldsymbol{q}$-Hypergeometric Orthogonal Polynomials\\ with $\boldsymbol{q=-1}$}

\Author{Luis VERDE-STAR}

\AuthorNameForHeading{L.~Verde-Star}

\Address{Department of Mathematics, Universidad Aut\'onoma Metropolitana,\\
Iztapalapa, Mexico City, Mexico}
\Email{\href{mailto:verde@xanum.uam.mx}{verde@xanum.uam.mx}}
\URLaddress{\url{http://mat.izt.uam.mx/mat/}}

\ArticleDates{Received October 31, 2024, in final form September 18, 2025; Published online October 02, 2025}

\Abstract{We obtain some properties of a class $\mathcal{A}$ of $q$-hypergeometric orthogonal polynomials with $q=-1$, described by a uniform parametrization of the recurrence coefficients. We construct a class $\mathcal{C}$ of complementary $-1$ polynomials by means of the Darboux transformation with a shift. We show that our classes contain the Bannai--Ito polynomials and their complementary polynomials and other known $-1$ polynomials. We introduce some new examples of~$-1$ polynomials and also obtain matrix realizations of the Bannai--Ito algebra.}

\Keywords{hypergeometric orthogonal polynomials; recurrence coefficients; $-1$ orthogonal polynomials; Bannai--Ito polynomials; Bannai--Ito algebra}

\Classification{33C45; 33D45}

\section{Introduction}
The classes of $q$-hypergeometric sequences of orthogonal polynomials with $q$ not equal to a~root of unity or with $q=1$ have been extensively studied for a long time \cite{Hyp}. In the case with $q=-1$, some polynomial sequences have been obtained as limits as $q$ goes to $-1$ of well-known $q$-orthogonal polynomials, such as the big and little $q$-Jacobi polynomials \cite{VZ1, VZ2}, or~as~transformations of the Bannai--Ito polynomials \cite{GenVZ2, GenVZ}. See also \cite{CohlCS, Tsuji, TsujiDunkl, TsuDualHahn}.

In \cite{Uni}, we showed that all the $q$-orthogonal polynomials sat\-isfy a gen\-er\-al\-ized dif\-fer\-ence-eigenvalue equation of order one, with respect to a Newtonian basis for the space ${\mathbb C}[t]$. In~\cite{qdiscr}, we~proved that all the sequences in the $q$-Askey scheme have a discrete orthogonality. Results about the discrete orthogonality of a subclass of the $-1$ polynomials will be presented elsewhere.

In \cite{Mops, Rec, PSGH,Uni}, we have used a matrix approach to study diverse aspects of the orthogonal polynomial sequences that produced general results that suggested new ways to classify the hypergeometric and basic hypergeometric orthogonal polynomial sequences. See \cite{KAsk, K2, K3}.

In the present paper, using some results from \cite{Uni} we construct a class $\cA$ of $q$-orthogonal polynomials with $q=-1$ and obtain a uniform parametrization for their recurrence coefficients. These polynomials may be considered as the elements of a $-1$-Askey scheme because they are constructed with the same procedure that we used to construct all the polynomial sequences in the Askey and the $q$-Askey schemes. The continuous part of the $-1$-Askey scheme was constructed in \cite{PellVZ}. We show that giving particular values for the parameters we can obtain some of the known $-1$ polynomials, such as the Bannai--Ito polynomials and the polynomials obtained as limits of the big and little Jacobi polynomials. We also obtain some examples of new $-1$ polynomials that seem interesting and deserve further study.

We also construct a class $\cC$ of $-1$ polynomials that we call the class of complementary~$-1$ polynomials because it contains the complementary Bannai--Ito polynomials, see \cite{GenVZ2,TsujiDunkl}. The elements of this class are Darboux transformations with a shift of the $-1$ polynomials in the class~$\cA$.
We describe a family of continuous $-1$ polynomials that contains the continuous Bannai--Ito polynomials obtained in \cite{PellVZ}. We also obtain some matrix realizations of the Bannai--Ito algebra.

In Section~\ref{section2}, we present a brief account of the construction of the bispectral $q$-hypergeometric orthogonal polynomials. In Section~\ref{section3}, we consider the case with $q=-1$ and present some examples. In Section~\ref{section4}, we introduce a parametrization that reduces the number of parameters and present some new examples. In Section~\ref{section5}, we consider a subclass of $\cA$ of polynomials that have simple recurrence coefficients. In Section~\ref{section6}, we construct the class $\cC$ of complementary~$-1$ polynomials and study a subclass of polynomial sequences whose recurrence coefficients are simple rational functions. In Section~\ref{section7}, we construct some matrix realizations of the Bannai--Ito algebra.
Finally, in Section~\ref{section8}, we mention some topics for further research.

\section[Construction of bispectral hypergeometric orthogonal polynomials]{Construction of bispectral hypergeometric\\ orthogonal polynomials}\label{section2}

In this section, we present a brief account of the construction of the
 bispectral hypergeometric orthogonal polynomials presented in our previous paper \cite{Uni}. Some similar results were obtained by Terwilliger in \cite{Terw,Ter3,Ter2} using a different approach.

Consider the linear difference equation
\begin{equation}\label{eq:diffeq}
 s_{k+3} = z ( s_{k+2} -s_{k+1}) + s_k, \qquad k\ge 0,
\end{equation}
where $z$ is a nonzero complex number and $s_k$ is a sequence of complex numbers with initial terms $s_0$, $s_1$, $s_2$. Since the characteristic polynomial is $t^3 - z t^2 + z t -1$, we see that the product of the roots is equal to one, the sum of the roots is equal to $z$, and 1 is a root. Therefore,
the characteristic roots of the difference equation are $1$, $q$, and $q^{-1}$, where $1 + q + q^{-1}=z$. If~${z=3}$, then $q=1$ is a triple root. If $z=-1$, then $q=-1$ is a double root, and if $z\ne 3$ and~${z\ne -1}$, then the roots $1$, $q$, $q^{-1}$ are distinct. Therefore, the general solution of \eqref{eq:diffeq} is of the form
$s_k= d_0 + d_1 q^k + d_2 q^{-k}$ when $z\ne 3$ and $z \ne -1$ and becomes $s_k=d_0 + d_1 k + d_2 k^2$ when $z=3$, and $s_k= d_0 +d_1 (-1)^k + d_2 k ( -1)^k$ when $z=-1$. These lattices were considered by Bochner and Hahn in their classification of the classical and the $q$-orthogonal polynomials. Such lattices also appear in~\cite{NSU}, where they are obtained from the solutions of a second order difference equation~\mbox{\cite[equation~(3.1.12)]{NSU}} that must satisfy certain additional conditions.

Let $x_k$, $h_k$, and $e_k$ be 3 solutions of \eqref{eq:diffeq}. These sequences will be used to construct a basis for the space ${\mathbb C}[t]$, a linear operator $\cD$ on the space of polynomials, and a sequence of orthogonal polynomials $u_k(t)$ that are eigenfunctions of $\cD$ with eigenvalues $h_k$. The sequence $e_k$ is used to simplify the definition of the operator $\cD$ and provides additional parameters that are needed to construct a theory that can be applied to all the hypergeometric and basic hypergeometric orthogonal polynomials in the Askey schemes.

The sequence $x_k$ determines the Newtonian basis $\{v_n(t) \mid n\ge 0\}$
of the complex vector space~${\mathbb C}[t]$ of polynomials in $t$, defined by
\begin{equation}\label{eq:Newton}
	v_n(t)= (t-x_0)(t-x_1)\cdots (t-x_{n-1}), \qquad n \ge 1,
\end{equation}
and $v_0(t)=1$.

We define the sequence $g_k$ by
\begin{equation*} 
g_k= x_{k-1} (h_k- h_0) + e_k, \qquad k \ge 1,
\end{equation*}
and $g_0=0$.
This sequence satisfies a linear difference equation of order five. We add the sequence $e_k$ to avoid some complicated restrictions on the initial values of $g_k$. Let us suppose that $h_k \ne h_j$ if $k \ne j$, and $g_k \ne 0$ for $k \ge 1$.

We use the basis $\{v_k\mid k \ge 0\}$ to define the linear operator $\cD$ by
\begin{equation*}
\cD v_k = h_k v_k + g_k v_{k-1}, \qquad k \ge 1.
\end{equation*}
	Since $g_0=0 $, we see that $\cD t^n$ is equal to $ h_n t^n $ plus a polynomial of lower degree.
The operator~$\cD$ is a generalized difference operator and can be realized
as a Dunkl type difference-reflection operator when $q=-1$, see~\cite{GenVZ2}.

	For $n \ge 0$, we define $u_n$ as the monic polynomial of degree $n$ which is an eigenfunction of $\cD$ with eigenvalue $h_n$. That is,
\begin{equation}\label{eq:eigenEq}
\cD u_n = h_n u_n, \qquad n \ge 0.
\end{equation}

 In \cite[p.~249]{Uni}, we showed that
\begin{equation*}
 u_n(t) = \sum_{k=0}^n c_{n,k} v_k(t), \qquad n \ge 0,
\end{equation*}
 where the coefficients $ c_{n,k}$ are given by
 \begin{equation}\label{eq:cnk}
c_{n,k}=\prod_{j=k}^{n-1} \frac{g_{j+1}}{h_n -h_j}, \qquad 0 \le k \le n-1,
 \end{equation}
and $c_{n,n}=1$ for $n \ge 0$.
This expression for $u_n(t)$ was also obtained in \cite{VZ} using a different approach. The idea of representing orthogonal polynomials in terms of a Newtonian basis was introduced by Geronimus in \cite{Ger}.

The matrix $C=[c_{n,k}]$, where the coefficients $c_{n,k}$ are defined in \eqref{eq:cnk}, is lower triangular and all its entries in the main diagonal are equal to $1$. Therefore, $C$ is invertible.  Let $ C^{-1}=[{\hat c}_{n,k}]$. Using some properties of divided differences, we proved in \cite[p.~251]{Uni} that
\begin{equation*}
  \hc_{n,k}=\prod_{j=k+1}^n \frac{g_j}{h_k - h_j}, \qquad 0 \le k \le n-1,
\end{equation*}
   and $\hc_{n,n}=1$ for $n \ge 0.$

	 The entries in the 0-th column of $C^{-1}$ are
\begin{equation}\label{eq:moments}
 \hc_{n,0} =\prod_{k=1}^n \frac{g_k}{h_0-h_k}, \qquad n \ge 1,
	\end{equation}
and $\hc_{0,0}=1$. We denote them by $m_n=\hc_{n,0}$ for $ n \ge 0$.

In \cite{Uni}, we also proved that the monic polynomial sequence $u_n(t)$ satisfies a three-term recurrence relation of the form
\begin{equation}\label{eq:3term}
	u_{n+1}(t)= (t-\beta_n) u_n(t) - \alpha_n u_{n-1}(t), \qquad n \ge 1,
\end{equation}
where the coefficients are given by
\begin{gather}\label{eq:beta}
	\beta_n= x_n + \frac{g_{n+1}}{h_n - h_{n+1}} -\frac{g_n}{h_{n-1} - h_n}, \qquad n\ge 0,
\\
\label{eq:alpha}
 \alpha_n = \frac{g_n}{h_{n-1} -h_n} \biggl(\frac{g_{n-1}}{h_{n-2} - h_n} - \frac{g_n}{h_{n-1} - h_n} + \frac{g_{n+1}}{h_{n-1}-h_{n+1}} +x_n - x_{n-1} \biggr), \qquad n\ge 1.\!\!
\end{gather}
Since $g_0=0$, the terms in the previous equations where $h_{-1}$ appears are equal to zero, see also~\mbox{\cite[Remark~2.5]{K3}}.

If all the $\alpha_n$ are positive and the $\beta_n$ are real, then the sequence $u_n$ is orthogonal with respect to a positive measure, and if all the $\alpha_n$ are nonzero, then $u_n$ is orthogonal with respect to a not necessarily positive definite moments functional.

The numbers $m_n$ are the generalized moments of the polynomial sequence $u_k(t)$ with respect to the Newtonian basis $\{ v_k(t)\mid k \ge 0 \}$ defined in \eqref{eq:Newton}.
From \eqref{eq:moments}, we see that $m_n$ satisfies the recurrence relation
\[
	 m_{n+1}=\dfrac{g_n}{h_0 -h_n} m_n, \qquad n\ge 1.
 \]

Let \smash{$\bigl\{p^{[1]}, p^{[2]}, p^{[3]}\bigr\}$} be the basis of solutions of \eqref{eq:diffeq} whose initial terms are $(1,0,0)$, $(0,1,0)$, and $(0,0,1)$, respectively.
These basic solutions are sequences of polynomials in $z$. It is easy to see that the initial terms of $p^{[1]}$ and $p^{[2]}$ are
\begin{gather*}
1,0,0,1,z,z (z-1), (z-1)(z^2-z-1),z (z-2)(z^2-z-1), \dots,
\\
0,1,0,-z, -(z-1)(z+1), -z(z^2-z-1), -z(z-1)(z-2)(z+1),\\
 -(z^2-z-1)(z^3-2 z^2 -z+1), \dots,
\end{gather*}
and that \smash{$p^{[3]}_k=p^{[1]}_{k+1}$} for $k \ge 0$.

Let us note that for $k \ge 0$ we have
\begin{gather*}
	x_k = x_0 p^{[1]}_k + x_1 p^{[2]}_k + x_2 p^{[3]}_k, \qquad
	h_k = h_0 p^{[1]}_k + h_1 p^{[2]}_k + h_2 p^{[3]}_k, \\
	e_k = e_0 p^{[1]}_k + e_1 p^{[2]}_k + e_2 p^{[3]}_k.
\end{gather*}
Therefore, $g_k$, $\alpha_k$, $\beta_k$, $c_{n,k}$ and ${\hat c}_{n,k}$ are functions of $z$ and the initial terms of the sequences~$x_k$,~$h_k$ and $e_k$.

The orthogonal polynomials with $q=1$ are obtained when we put $z=3$.
Taking $z=1+ q+q^{-1}$, where $q\ne \pm 1$, we obtain all orthogonal polynomials in the $q$-Askey scheme. In the following section we consider the case $q=-1$.

The properties of the polynomial sequence represented by the matrix $C$ can be expressed in terms of matrix equations as we describe next.
Let $X$ be the shift matrix defined by $X_{k,k+1}=1$ for $k\ge 0$ and all the other entries equal zero.
 The matrix $S$ is the transpose of $X$. The diagonal matrix $H=\diag(h_0,h_1,h_2,\dots)$ is called the matrix of eigenvalues. The diagonal matrix $F=\diag(x_0,x_1,x_2,\dots)$ is the matrix of nodes associated with the Newtonian basis. The diagonal matrix $G=\diag(g_1,g_2,g_3,\dots)$ is associated with the operator $\cD$. The Jacobi matrix~$L$ is defined by $L_{k,k+1}=1$, $L_{k,k}=\beta_k$, $L_{k+1,k}=\alpha_{k+1}$, for $k\ge 0$.

 With respect to the Newtonian basis the operator of multiplication by the independent variable $t$ is represented by $ X+F$, and the operator $\cD$ by $( H + S G)$.
 The three-term recurrence relation \eqref{eq:3term} corresponds to the equation
 \begin{equation*}
	 LC=C(X+F)
 \end{equation*}
and the difference-eigenvalue equation \eqref{eq:eigenEq} is expressed in terms of matrices as
\begin{equation*}
	C(H+SG)=HC.
\end{equation*}
See \cite{Uni} for a more detailed account of these matrix equations.

\section[The polynomials with q=-1]{The polynomials with $\boldsymbol{q=-1}$}\label{section3}

We consider now the class $\cA$ of orthogonal polynomials obtained when the parameter $z$ in the difference equation \eqref{eq:diffeq} is equal to $-1$, and thus the roots of the characteristic polynomial of the difference equation are 1 and $-1$, which is a double root.

When $z=-1$ the basic solutions of \eqref{eq:diffeq} become the sequences
\begin{gather*}
1,0,0,1,-1,2,-2,3,-3,4,-4,\dots, \\
0,1,0,1,0,1,0,1,0,1,0, \dots, \qquad \hbox{and} \\
0,0,1,-1,2,-2,3, -3, 4,-4,5, \dots,
\end{gather*}
and a solution with initial values $s_0$, $s_1$, $s_2$ has the form
\[
s_0, s_1, s_2, s_0+s_1-s_2, -s_0+2 s_2, 2 s_0 + s_1 - 2 s_2, \dots.
\]
Another basis for the space of solutions of the difference equation \eqref{eq:diffeq} when $z=-1$ is $\bigl\{ 1, (-1)^k, (-1)^k k\bigr\}$. We write the sequences $x_k$, $h_k$ and $e_k$ in terms of this basis as follows:
\begin{gather}
 x_k = b_0 + b_1 (-1)^k + b_2 (-1)^k k, \\
 h_k = a_0 + a_1 (-1)^k + a_2 (-1)^k k, \\
 e_k = d_0+ d_1 (-1)^k + d_2 (-1)^k k. \label{eq:newparam}
\end{gather}
It is clear that the coefficients in these representations can be written in terms of the corresponding initial terms, and vice-versa.
We will see that the number of parameters needed to describe the orthogonal polynomials with $q=-1$ can be reduced to four.

The coefficient $b_0$ is the constant part of the sequence of nodes $x_k$ and hence, a change in the value of $b_0$ is equivalent to a translation of the polynomials $v_k(t)$ in the Newtonian basis and also of the orthogonal polynomials $u_k(t)$. Since translations do not change the relevant properties of the polynomials $u_k(t)$, sometimes we take $b_0=0$ in order to simplify the expressions for the coefficients of the $u_k(t)$ and of their three-term recurrence relation.

Let us note that the eigenvalues $h_k$ appear only in terms of the form $h_k-h_n$ in the formulas of the previous section. Therefore,
the coefficient $a_0$, which is the constant part of the eigenvalues~$h_k$, is always cancelled and thus we take $a_0=0$.

In order to have $g_0=0$, we need $e_0=0$ and thus we put $d_0= -d_1$ in \eqref{eq:newparam}.
If we take $b_0=0$ ,then there remain six parameters, $a_1$, $a_2$, $b_1$, $b_2$, $d_1$, $d_2$.

The sequence $g_k=x_{k-1}(h_k - h_0)+ e_k$ is in this case
\begin{gather*}
	g_k= -k ( (k-1)a_2 b_2 -d_2+ a_2 b_1) \qquad \hbox{if $k$ is even},\\
	g_k= a_2 b_2 k^2 + (d_2+ (b_1-b_2) a_2 + 2 a_1 b_2) k +2 ( b_1- b_2) a_1 + 2 d_2\qquad \hbox{if $k$ is odd}.
\end{gather*}

For even $n$, the recurrence coefficients are
\begin{gather}\label{eq:alphaev}
	\alpha_n= \dfrac{n ((n-1)a_2 + 2 a_1) (n a_2 b_2+2 a_1 b_2 -a_2 b_1 + d_2)((n-1) a_2 b_2 + a_2 b_1 -d_2)}{a_2( (2 n-1) a_2 + 2 a_1)^2 },
\\
\label{eq:betaev}
\beta_n= \dfrac{a_2(2 a_1 b_2 +a_2(b_2-2 b_1)-4 d_1) n - (2 d_1+d_2)( 2 a_1 -a_2)}{((2n-1) a_2 + 2 a_1) ((2n+1) a_2 + 2 a_1)}.
\end{gather}

For odd $n$, we have
\begin{equation}\label{eq:alphaodd}
	\alpha_n= \dfrac{-w_1(n) w_2(n) }{a_2( (2 n-1) a_2 + 2 a_1)^2 },
\end{equation}
where
\begin{gather*}
w_1(n)=a_2 b_2 n^2+ ((b_1-b_2) a_2 + 2 a_1 b_2 + d_2) n + 2 (b_1-b_2) a_1 + 2 d_1,\\
w_2(n)=  a_2^2 b_2 n^2 + \bigl(- (b_1+ b_2) a_2^2 + ( 2 a_1 b_2 - d_2) a_2\bigr) n + a_2^2 b_1\\
\hphantom{w_2(n)=}{} + (2 d_1+ d_2- 2 a_1 b_2) a_2 - 2 a_1 d_2,
\end{gather*}
and
\begin{equation}\label{eq:betaodd}
	\beta_n= \dfrac{- w_3(n)}{((2 n-1) a_2 + 2 a_1) ( ( 2 n +1) a_2 + 2 a_1) },
\end{equation}
where
\begin{align*}
w_3(n)={}& a_2( 2 a_1 b_2 + a_2(b_2- 2 b_1)-4 d_1) n - 4 a_1^2 b_2\\
&{+}\,((4 b_1 -2 b_2) a_2 + 4 d_1- 2 d_2) a_1 + ( 2 d_1 + d_2) a_2.
\end{align*}

Giving appropriate values to the parameters in the equations \eqref{eq:alphaev}--\eqref{eq:betaodd} we can obtain the recurrence coefficients of the $-1$ orthogonal polynomial sequences in the class $\cA$. For example, the recurrence coefficients for the Bannai--Ito polynomials \cite{GenVZ} are obtained with
\begin{gather*}
a_2=1,\qquad b_2=1, \qquad a_1=\alpha+\beta +\gamma+ \delta+ 3/2 ,\qquad b_1=1 + 2 \beta, \\
d_1= 2 \gamma \delta -2 \alpha \beta - \beta + \gamma + \delta +1/2 , \qquad d_2=-1- 2 \alpha.
\end{gather*}
Here $\alpha$, $\beta$, $\gamma$, $\delta$ are the parameters used in \cite{GenVZ} to define the Bannai--Ito polynomials.

If we take
\begin{gather*}
b_1=c,\qquad b_2=0,\qquad a_1=\dfrac{1 +\alpha + \beta}{2 ( 1+c)},\qquad a_2=\dfrac{1}{1+c},\\
 d_1=\dfrac{1 +\alpha + \beta}{2 ( 1+c)}-\dfrac{1+\alpha}{2},\qquad d_2= \dfrac{1}{1+c},
 \end{gather*}
we get the recurrence coefficients of the polynomial sequence obtained from the big $q$-Jacobi polynomials by taking the limit as $q$ goes to $-1$. See \cite{CohlCS} and \cite{VZ2}.

For the $-1$ polynomials obtained in \cite{VZ1} as limits of the little $q$-Jacobi polynomials, we get the recurrence coefficients with
\[
b_1=0,\qquad b_2=0,\qquad a_1=\dfrac{1+\alpha + \beta}{2},\qquad a_2=1,\qquad d_1=\dfrac{\alpha}{2},\qquad d_2=1.
\]
Note that these coefficients are obtained from those of the big $q$-Jacobi polynomials taking $c=0$ and $\beta=\alpha$.

We obtain $-1$ polynomials related with Chebyshev polynomials as follows. Taking
\[
a_1=\dfrac{a_2}{2},\qquad b_1=0,\qquad b_2=0,\qquad d_1=0, \qquad d_2=a_2,
\]
the recurrence coefficients are $\alpha_n=1/4$ for $n\ge 1$, $\beta_0=-1/2$, and $\beta_n=0$ for $n\ge 1$. The corresponding monic orthogonal polynomial sequence $u_n(t)$ is related with the monic Chebyshev polynomials of the first kind $\tilde T_n(t)$ by
\[
u_n(t)= \sum_{k=0}^n \dfrac{\tilde T_k(t)}{2^{n-k}}, \qquad n \ge 0.
\]
If instead of $d_2= a_2$, we put $d_2=- a_2$ then we get $\beta_0=1/2$ and the corresponding orthogonal polynomials $w_n(t)$ satisfy
\[
w_n(t)=\sum_{k=0}^n (-1)^{n-k} \dfrac{\tilde T_k(t)}{2^{n-k}}, \qquad n \ge 0.
\]

Another simple case is obtained when we take
\begin{gather*}
 a_1=1/2,\qquad a_2=1,\qquad b_1=1/2,\qquad b_2=1,\qquad d_1=1/4,\qquad d_2=-1/2.
 \end{gather*}
In this case, we have $\alpha_n=-n^2/4$, for $ n\ge 1$, and $\beta_n=0$, for $n\ge 0$. These examples provide matrices that can be used to obtain generators of the Bannai--Ito algebra.

We can also obtain the recurrence coefficients of other $-1$ polynomials considered in \cite{PellVZ} by taking appropriate values for the parameters $a_1$, $a_2$, $b_1$, $b_2$, $d_1$, $d_2$, without taking limits.

A classification of the polynomial sequences in the class $\cA$ can be obtained by the method used by Koornwinder to construct the schemes in \cite{KAsk, K2,K3}.

\section[A change of parameters for the case with b\_2 ne 0]{A change of parameters for the case with $\boldsymbol{b_2 \ne 0}$}\label{section4}
In this section, we consider the class of orthogonal $-1$ polynomials for which the parameter $b_2$ is nonzero. Since $x_k=b_1 (-1)^k + b_2 k (-1)^k$ we can see that when $b_2$ is nonzero the nodes $x_k$ are pairwise distinct.
 We will introduce a parametrization for the recurrence coefficients that simplifies the expressions and also reduces the number of independent parameters.

From equations \eqref{eq:alphaev} and \eqref{eq:alphaodd}, we see that $a_2$ must be nonzero. This condition is also needed to have $h_k \ne h_n$ if $ k \ne n$.

We define the parameters $r$, $s$, $t_1$ and $t_2$ by the equations
\begin{gather*}
a_1=\biggl(s + \dfrac{1}{2}\biggr) a_2,\qquad b_1=\biggl(r+\dfrac{1}{2}\biggr) b_2,\qquad d_1=\biggl(\dfrac{t_1}{2} +\dfrac{1}{4}\biggr) a_2 b_2, \qquad d_2=\biggl( t_2- \dfrac{1}{2}\biggr) a_2 b_2.
\end{gather*}
Substitution of $a_1$, $b_1$, $d_1$, $d_2$ in equations \eqref{eq:alphaev}--\eqref{eq:betaodd} yields, if $n$ is even
\begin{gather}
\alpha_n = \dfrac{-b_2^2n ( n+ 2 s) ( n+ 2 s -r + t_2) ( n+ r - t_2)}{4 (n+s)^2},\\
	\beta_n =\dfrac{-b_2((r-s+t_1) n +s (t_1+t_2))}{2 (n+s)(n+1+s)},\label{eq:alphaevrs}
\end{gather}
 and if $n$ is odd
\begin{align}
\alpha_n ={}& \dfrac{-b_2^2\bigl(n^2+(r+2s+t_2) n + (2 r -1)s + r + t_1\bigr)}{4(n+s)^2} \\
 	 &{}{\times} \, \bigl(n^2+(2s -r - t_2) n -(2 t_2 + 1)s + r + t_1\bigr),\label{eq:alphaoddrs}
\end{align}
and
\begin{equation}\label{eq:betaoddrs}
	\beta_n= \dfrac{b_2\bigl((r-s+t_1) n -2 s^2 +( 2r +t_1 - t_2-1)s + r + t_1\bigr)}{2(n+s)( n+1+s)}.
\end{equation}

Let us note that $a_2$ does not appear in the new expressions for the recurrence coefficients. Note also that, for all $n$, $b_2$ is a factor of $\beta_n$ and $b_2^2$ is a factor of $\alpha_n$. In \cite[Appendix, p.~215]{Chi}, Chihara shows that the effect of such factors on the corresponding monic polynomial sequences is a change of the independent variable of the form $w_n(x)=b_2^n u_n(x/b_2)$. If we consider that two monic polynomial sequences related in such way are equivalent, then $b_2$ can be taken as any nonzero number. Therefore, the recurrence coefficients are essentially determined by the four parameters $r$, $s$, $t_1$, $t_2$, see also \cite[Remark~2.5]{K3}.

We describe next the recurrence coefficients of a family of continuous $-1$ polynomial sequences. That is, such that all $\alpha_n$ are positive and all $\beta_n$ are real.

Since the leading term in the numerators of \eqref{eq:alphaevrs} and \eqref{eq:alphaoddrs} is $- b_2^2 n^4$, it is clear that $b_2^2$ must be a negative number if we want $\alpha_n >0$ for $n\ge 1$. We take $b_2={\rm i}$, where ${\rm i}^2=-1$.
Let $y=y_1 + {\rm i} y_2$ and $w=w_1 + {\rm i} w_2$ be complex numbers.

Substitution of{\samepage
\[ b_2={\rm i},\qquad s=y_1,\qquad r=\dfrac{y-w}{2},\qquad t_2=-\dfrac{y+w}{2},\qquad t_1=\dfrac{\bar y +w}{2} -(y_1 y_2 - w_1 w_2) {\rm i}, \]
 in \eqref{eq:alphaevrs}--\eqref{eq:betaoddrs} yields the following formulas.}

 If $ n$ is even,
 \begin{gather}\label{eq:alphaevc}
\alpha_n=\dfrac{n ( n+2 y_1) ( n + y) ( n+\bar y)}{4 (n+ y_1)^2},
\\ \label{eq:betaevc}
	 \beta_n=\dfrac{( w_1 w_2 -y_1 y_2) n +y_1 ( w_1 w_2 - y_1 y_2) - y_1 y_2}{2 (n + y_1) ( n+1 + y_1)}.
 \end{gather}
If $n$ is odd,
\begin{gather}\label{eq:alphaoddc}
	\alpha_n=\dfrac{ (n+y_1+w_1) (n+y_1 - w_1) (n+ y_1 + {\rm i} w_2) ( n+ y_1 - {\rm i} w_2)}{4 ( n+y_1)^2},
\\ \label{eq:betaoddc}
	 \beta_n=\dfrac{(y_1 y_2 - w_1 w_2) n + y_1( y_1 y_2 - w_1 w_2) -w_1 w_2}{2 ( n+y_1) ( n+1 + y_1)}.
\end{gather}

It is easy to see that $\alpha_n >0$ for all $n \ge 1$ if $ y_1>-1$ and $ w_1 + y_1 > -1$.

Let us consider some simple examples. If we take $y_1=0$ and $w_1=0$, then for even $n$ we have \smash{$\alpha_n= \bigl( n^2 + y_2^2\bigr)/4$}, and for odd $n$ we have \smash{$\alpha_n= \bigl( n^2 + w_2^2\bigr)/4$}, $\beta_0=-y_2/2$ and $\beta_n=0$ for $n\ge 1$.

If $y=w=0$, then $ \alpha_n=n^2/4$ for $n\ge 1$ and $\beta_n=0$ for $n\ge 0$.

If $y_2=0$ and $w_2=0$, then for even $n$ we have $\alpha_n=n ( n+ 2 y_1)/4$, for odd $n$ we have $\alpha_n=(n+y_1 - w_1) (n+ y_1 - w_1)/4$, and $\beta_n=0$ for $n\ge 0$.

The recurrence coefficients of the continuous Bannai--Ito polynomials, which are described in~\cite{PellVZ} using parameters $\alpha$, $\beta$, $\gamma$, $\delta$, are obtained by substitution of
\[ y_1=1+ 2 ( \alpha + \gamma),\qquad y_2=2 ( \beta + \delta),\qquad w_1=2 ( \gamma - \alpha), \qquad w_2=2 ( \delta - \beta), \]
 in \eqref{eq:alphaevc}--\eqref{eq:betaoddc}.

 There are also continuous $-1$ polynomial sequences in the class with $b_2=0$.
 We give some examples in Section~\ref{section5}.

 We obtain next another useful parametrization of the class of $-1$ polynomials with $b_2 \ne 0$ that gives factorized expressions for all the $\alpha_n$.
 Let $p$ and $q$ be roots of the numerator in \eqref{eq:alphaoddrs} and substitute
\[ t_1=p q +s -r ( 1+ 2s), \qquad t_2=-(p+q+r+2s), \]
 in \eqref{eq:alphaevrs}--\eqref{eq:betaoddrs}.
 We obtain, for even $n$,
\begin{gather*}
\alpha_n=-b_2^2\dfrac{n ( n+ 2s) (n-p-q- 2 r) ( n+p+q+2r+ 2s)}{4(n+s)^2},
\\ 
\beta_n=- b_2\dfrac{(pq-2 rs) n - s^2 ( 1+2r) +s (pq-p-q-2r)}{2 (n+s) ( n+1+s)},
 \end{gather*}
 and for odd $n$
\begin{gather*}
	\alpha_n=-b_2^2\dfrac{(n-p) (n-q) ( n+p+ 2s) (n+q + 2 s) }{4(n+s)^2},
\\ 
\beta_n= b_2\dfrac{(pq-2 rs) n + p q+ s^2 ( 1-2r) +s (pq+p+q)}{2 (n+s) ( n+1+s)}.
 \end{gather*}

 \section[A parametrization for the -1 polynomials with b\_2=0]{A parametrization for the $\boldsymbol{-1}$ polynomials with $\boldsymbol{b_2=0}$}\label{section5}

 If $b_2=0$ and $b_1 \ne 0$, then $x_k= b_1 (-1)^k$, for $k\ge 0$, and hence there are only two distinct nodes. If both, $b_1$ and $b_2$ are zero then $x_k=0$ for $k\ge 0$. Recall that $a_2$ must be nonzero.

 For the case with $b_2=0$, we obtain a parametrization of the recurrence coefficients by substitution of
 \[ b_2=0,\qquad a_1=(s+ 1/2) a_2, \qquad d_1=t_1 a_2/2, \qquad d_2= t_2 a_2, \]
 in \eqref{eq:alphaev}--\eqref{eq:betaodd}. We obtain, for even $n$,
\begin{gather}\label{eq:alphaevb2}
\alpha_n=\dfrac{(b_1-t_2)^2 n ( n+2 s)}{4 (n+s)^2},
\\
\label{eq:betaevb2}
\beta_n=-\dfrac{(b_1+ t_1) n + s ( t_1+ t_2)}{2(n+s) ( n+1+s)},
\end{gather}
and for odd $n$ we have
\begin{gather}\label{eq:alphaoddb2}
	\alpha_n=\dfrac{((b_1+t_2) n + 2 b_1 s+ b_1+t_1)( (b_1+ t_2)n+2 s t_2 - b_1 - t_1) }{4 (n+s)^2},
\\\label{eq:betaoddb2}
	\beta_n=\dfrac{ (b_1+ t_1) (n+1) + s(2 b_1+t_1-t_2) }{2(n+s) ( n+1+s)}.
\end{gather}

Let us note that the recurrence coefficients are expressed in terms of the four parameters~$b_1$,~$s$,~$t_1$ and~$t_2$.

We present next some examples obtained by giving particular values to the parameters in~\eqref{eq:alphaevb2}--\eqref{eq:betaoddb2}.
These examples are constructed so that the recurrence coefficients have simple expressions as functions of the index $n$. We also try to show that it is easy to obtain recurrence coefficients that are rational functions of the index $n$ with several different values for the degrees of the numerators and the denominators. These degrees may be used to obtain a classification criterion. The examples also show that we can construct families of $-1$ polynomials that are determined by several parameters.

Let us recall that the values of the parameters also determine the generalized difference operator $\cD$ and the Newtonian basis. With an appropriate change of bases the operator $\cD$ can be transformed to an operator with respect to the standard basis of monomials.

Additional examples are presented in Section~\ref{section7}, where we consider a subclass of the class of complementary $-1$ polynomials.

\begin{Example}
Taking $b_1=0$, $t_1=2s$ and $t_2=1$, we get
\begin{gather*}
\alpha_n=\dfrac{n(n+2s)}{4 (n+s)^2}, \qquad n\ge 1,
\qquad
	\beta_n=(-1)^{n-1}\dfrac{s(2 n + 2s+1)}{2 (n+s)(n+1+s)}, \qquad n\ge 0.
\end{gather*}
If we now put $s=0$ in these equations, we get $\alpha_n=1/4$ for $n\ge 1$, $\beta_0=-1/2$, and $\beta_n=0$ for $n\ge 1$. This is one of the examples related with the Chebyshev polynomials presented in~Section~\ref{section3}.
\end{Example}
\begin{Example}
With $s=0$ and $t_2=-b_1$, we obtain, for even $n$,
\[ \alpha_n=b_1^2, \qquad \hbox{and }\qquad \beta_n= -\dfrac{b_1+t_1}{2 (n+1)}, \]
and for odd $n$,
\[ \alpha_n=-\dfrac{(b_1+t_1)^2}{4 n^2}, \qquad \hbox{and } \qquad \beta_n=\dfrac{b_1+t_1}{2 n}. \]
\end{Example}
\begin{Example}
Here we take $b_1=0$ and $t_2=1$ and we get, for even $n$,
\[ \alpha_n=\dfrac{n ( n+2s)}{4(n+s)^2}, \qquad \hbox{and} \qquad \beta_n=-\dfrac{2 t_1 n + s(1+ 2 t_1)}{2(n+s)(n+1+s)}, \]
and for odd $n$
\[ \alpha_n=\dfrac{(n+2s-2t_1)(n+ 2 t_1)}{4(n+s)^2}, \qquad \hbox{and} \qquad \beta_n=
\dfrac{2 t_1( n+1) + s( 2 t_1-1)}{2(n+s)(n+1+s)}. \]
If we now take $s=0$ and $t_1=1/2$, we obtain, for even $n$,
\[ \alpha_n=\frac{1}{4}, \qquad \hbox{and} \qquad \beta_n=-\dfrac{1}{2 (n+1)}, \]
and for odd $n$
\[ \alpha_n=\dfrac{n^2-1}{4 n^2}, \qquad \hbox{and} \qquad \beta_n=\dfrac{1}{2 n}. \]
\end{Example}

\section[The class of complementary -1 polynomials]{The class of complementary $\boldsymbol{-1}$ polynomials}\label{section6}

In this section, we describe a class of orthogonal polynomial sequences whose monic Jacobi matrices are obtained as Darboux transformations with a shift $w$ of the Jacobi matrices of sequences in the class $\cA$ of $-1$ polynomials. We consider first Darboux transformations for $q$-polynomial sequences with $q \ne -1$.

Let $w$ be a complex number and let $L$ be the Jacobi matrix defined in Section~\ref{section2}. Let ${Y=\diag(y_0,y_1,y_2,\dots)}$ and $Z=\diag(z_0,z_1,z_2,\dots)$ be the unique diagonal matrices that satisfy
\begin{gather} \label{eq:LUfactor}
L -w I = (I + Y S) ( X+Z).
\end{gather}
From this factorization, it is easy to see that
\begin{equation*}
	L_{k+1,k}=\alpha_{k+1}= y_{k+1} z_k, \qquad k\ge 0,
\qquad 
	L_{k,k}=\beta_k= y_k+ z_k +w, \qquad k\ge 0,
\end{equation*}
where the entries of $L$ are given in \eqref{eq:alpha} and \eqref{eq:beta}.

The entries $y_k$ and $z_k$ of the matrices $Y$ and $Z$
 can be obtained from \eqref{eq:LUfactor}. They are rational functions of $q^k$, $w$ and the parameters in $\alpha_k$. As functions of $q^k$ the degrees of the numerator and the denominator of $y_k$ and $z_k$ are increasing functions of $k$, but the degrees of the numerator and the denominator of $\alpha_k$ as functions of $q^k$ are independent of $k$, and they are both equal to eight. For $\beta_k$, the corresponding degrees are both equal to four, see \cite[Section~7]{Uni}. This happens because, for each $k$, $y_{k+1}$ has a rational factor whose reciprocal is a factor of $z_k$, and thus such factors cancel in the product $y_{k+1} z_k$.

 We define the matrix
 \begin{equation}\label{eq:Darb}
	 M= (X+Z) (I+ Y S)+ w I.
	\end{equation}
This matrix is a monic Jacobi matrix and it is called the Darboux transform of $L$ with shift $w$, see \cite{BuenoM}.

From the definition of $M$, we can see that
\begin{gather}\label{eq:alMyz}
	M_{k+1,k}= y_{k+1} z_{k+1}, \qquad k \ge 0,\\
\label{eq:beMyz}
	M_{k,k}= y_{k+1} + z_k + w, \qquad k\ge 0.
\end{gather}

For general values of the shift $w$, the cancellation of factors that occurs in the products $y_{k+1} z_k$ does not occur in the products $ y_{k+1} z_{k+1}$. Because of this fact, the entries of $M$ are rational functions of $q^k$ with numerators and denominators whose degrees are increasing functions of~$k$. Therefore, in such cases, the $q$-orthogonal polynomial sequence determined by the matrix~$M$ is not in the class $\cH_q$, defined in \cite[Section~7]{Uni}, which contains the sequences in the $q$-Askey scheme. For the particular value $w=x_0$, the Jacobi matrix $M$ determines a polynomial sequence in~$\cH_q$, and $M$ is obtained from $L$ by a simple modification of the parameters, that we will describe next.\looseness=-1

In the general case of the $q$-hypergeometric orthogonal polynomials, with $q\ne 1$ and ${q \ne -1}$, the characteristic roots of the difference equation \eqref{eq:diffeq} are $1$, $q$ and $q^{-1}$ and are distinct.
Therefore, the sequences $h_k$, $ x_k$, and $e_k$ are given by
\begin{gather*} 
h_k= a_0 + a_1 q^k + a_2 q^{-k}, \qquad
	x_k= b_0 + b_1 q^k + b_2 q^{-k}, \qquad
	e_k= d_0 + d_1 q^k + d_2 q^{-k},
\end{gather*}
where $a_0$ may be taken as zero and $d_0=-d_1 -d_2$.

The recurrence coefficients $\alpha_k$ and $\beta_k$ depend on the parameters $a_1$, $a_2$, $b_0$, $b_1$, $b_2$, $d_1$, $d_2$ and~$q$, and the entries of $M$ depend on the same parameters and also on $w$.

If we take $w=x_0=b_0+b_1+b_2$, then applying the substitution
\begin{equation}\label{eq:subpar}
	( a_1, a_2, b_0, b_1, b_2, d_1, d_2)\rightarrow \bigl( q a_1, a_2, q b_0, q^2 b_1, b_2, q^2 d_1, q d_2\bigr),
\end{equation}
to $(\alpha_k, \beta_k)$ we obtain $\bigl( q^2 M_{k+1,k}, q M_{k,k}\bigr)$, which corresponds to a re-scaling of the variable in the polynomial sequence determined by $M$. The value $w=x_0$ is the only one for which this result holds for arbitrary values of the parameters $a_1$, $a_2$, $b_0$, $b_1$, $b_2$, $d_1$, $d_2$. Note that \eqref{eq:subpar} is invertible and therefore the Darboux transformation with $w=x_0$ sends the class $\cH_q$ onto itself.

For the Askey--Wilson polynomials, due to the symmetries of the four parameters~$a$,~$b$,~$c$,~$d$ that are used to describe the recurrence coefficients, the initial node $x_0$ can be written as $x_0=\bigl(r+ r^{-1}\bigr)/2$, where $r$ is any of $a$, $b$, $c$, $d$, see \cite[equation~(14.1.5)]{Hyp} and \cite[equation~(7.15)]{Uni}.

For the class $\cH_1$ of $q$-orthogonal polynomials with $q=1$, defined in \cite[Section~8]{Uni}, the sequences $h_k$, $x_k$, $e_k$ are given by
\begin{gather*} 
h_k= a_1 k+ a_2 k^2, \qquad
x_k= b_0 + b_1 k + b_2 k^2, \qquad
e_k= d_1 k + d_2 k^2.
\end{gather*}
In this case, the Darboux transformation with shift $w$ of the generic Jacobi matrix $L$ yields in the general case a complicated Jacobi matrix $M$, but if we take $w=x_0=b_0$, then the substitution
\begin{equation}\label{eq:subpar1}
	( a_1,a_2, b_0, b_1, b_2, d_1, d_2)\rightarrow ( a_1+a_2, a_2, b_0+b_1+b_2, b_1+ 2 b_2, b_2, d_1+d_2, d_2),
\end{equation}
applied to $(\alpha_k, \beta_k)$ gives us $(M_{k+1,k}, M_{k,k})$.
Therefore, the Darboux transformation with shift $w=x_0$ sends $\cH_1$ onto itself, since the change of parameters \eqref{eq:subpar1} is invertible.

We consider next the Darboux transformation of the class $\cA$ of $-1$ polynomials. In this case, the recurrence coefficients are given in equations \eqref{eq:alphaevb2}--\eqref{eq:betaoddb2}.

Let us recall the main properties of the matrix $C$ of $-1$ orthogonal polynomials. We have the recurrence relation
\begin{equation}\label{eq:LC}
	LC=C(X+F),
\end{equation}
and the eigenvalue equation
\begin{equation}\label{eq:ChSg}
	C(H+ S G)=H C.
\end{equation}

Note that the matrix $I+YS$, defined in \eqref{eq:LUfactor} is invertible. Then we define
\begin{equation*}
	\tilde C= (I+ YS)^{-1} C.
\end{equation*}
From equation \eqref{eq:LUfactor}, we obtain
\begin{equation*}
	(I +YS)^{-1} L (I+YS)= (X+S) ( I+YS) + wI=M,
\end{equation*}
and from \eqref{eq:LC} we get
\begin{equation*}
	(I+YS)^{-1} L (I+YS) \tilde C = \tilde C ( X+F).
\end{equation*}
Therefore, we have
$
M \tilde C= \tilde C (X+F)$,
and hence $\tilde C$ is the matrix of orthogonal polynomials associated with $M$, expressed in terms of the Newtonian basis $\{v_n(t)\mid n \ge 0\}$.
Therefore, the polynomial sequences described by~$\tilde C$ are~$-1$ orthogonal polynomials. The set of all such polynomial sequences is the class of complementary~$-1$ orthogonal polynomial sequences that we denote by $\cC$.

Define the matrix
$
	B={\tilde C}^{-1} C (H+S G)C^{-1} \tilde C$.
Then by \eqref{eq:ChSg} we get
$
	\tilde C B = H \tilde C$.
The matrix $B$ is similar to the bidiagonal matrix $H+ S G$, but in the general case it is not a~banded matrix. The matrix $C^{-1} \tilde C$ represents a change of bases on the space of polynomials. For some particular cases, it may be possible to find matrices of change of bases that convert $B$ into a banded matrix that represents some kind of difference operator.

The entries of $M$ are rational functions of the parameters $a_1$, $a_2$, $b_0$, $b_1$, $b_2$, $d_1$, $d_2$, $w$ and the index $k$. In the general case the sequence of degrees, as functions of $k$ or $a_2$, of both the numerator and the denominator of $M_{k+1,k} $ for $k\ge 1 $ is $6,10,10,14,14,18,18,\dots $ and the analogous sequence for $M_{k,k}$ is $ 4,6,8,10,12,\dots$. This shows that the coefficients of the three-term recurrence coefficients satisfied by the complementary $-1$ orthogonal polynomial sequences are in general quite complicated.

Giving appropriate values to the parameters we can obtain matrices $M$ with entries whose numerators and denominators have bounded or constant degrees. For example, the values ${b_0=0}$, $b_1=-w$, $b_2=0$, $d_1=0$, $d_2=w a_2$, $a_1=(r/2 -1)a_2$ give us a matrix $M$ for which the denominator of $M_{k+1,k}$ is equal to $(r+k)(r+k+1)^2 (r+k+2)$, which has degree four for every~$k$, and the numerator also has degree four for every $k$.

There is a subclass of $\cC$ whose elements have simple recurrence coefficients, similar to those of the complementary Bannai--Ito polynomials.
This subclass is found by looking for relations among $w$ and some of the other parameters that yield cancellation of certain factors in the numerator and the denominator of $M_{k+1,k}$. We found that such cancellation occurs if ${w=b_0+b_1=x_0}$ or $w=-d_2/a_2$.

\subsection[The subclass C\_0]{The subclass $\boldsymbol{\cC_0}$}
Let us denote by $\cC_0$ the subclass of $\cC$ of all the sequences obtained with $w=b_0+b_1$.
When $w=b_0+b_1$ the sequences $y_k$ and $z_k$ become quite simple. They are rational functions with different formulas for even and odd values of $k$ that we present next.
 We define the following functions:
 \begin{gather*}
	 y_e(k)=\dfrac{2 k ( ( b_0-b_1+ 2 k b_2) a_2 + 2 a_1 b_2 + d_2)}{2 a_1 + (4 k -1) a_2},
\\ 
	y_o(k)=\dfrac{2 a_2^2 ( (2 k+1) b_2 -b_0 -b_1) k +a_2 (2 k a_1 b_2 -k d_2 + d_1) - a_1 d_2 }{ a_2 (2 a_1+ (4k+1) a_2) },
\\ 
	z_e(k)=- \dfrac{(2 k b_2 +b_0+b_1)(a_2 ( 2 k+1) + 2 a_1) + ( 2 k+1) d_2 + 2 d_1}{2 a_1 + ( 4 k+1) a_2},
\\ 
	z_o(k)=\dfrac{ (2 a_1+(2k+1) a_2)( a_2 ( b_0-b_1 - ( 2k+1)b_2) + d_2)}{a_2 ( 2 a_1+ ( 4k +3) a_2)}.
\end{gather*}
Then the sequences $y_k$ an $z_k$ for the subclass $\cC_0$ are given by
\begin{gather*}
 y_{2k}= y_e(k), \qquad  k\ge 0,\qquad
 y_{2k+1}=y_o(k),\qquad  k\ge 0,  \\
 z_{2k}=z_e(k),\qquad  k\ge 0,\qquad
 z_{2 k+1}=z_o(k), \qquad  k \ge 0.
\end{gather*}
Let us note that, for every $k \ge 0$, $y_k$ and $z_k$ are rational functions of $k$ with numerator of degree two and denominator of degree one.
Therefore, since $M_{k+1,k}= y_{k+1} z_{k+1}$ for every $k\ge 0$, the numerator of $M_{k+1,k}$ has degree four and its denominator has degree two.

Explicit formulas for $M_{k+1,k}$ for the subclass $\cC_0$ are obtained immediately from
\begin{gather*}
 M_{2k+1,2k}= y_o(k) z_o(k),\qquad  k\ge 0,\qquad
 M_{2k, 2k-1}= y_e(k) z_e(k), \qquad  k\ge 1.
\end{gather*}
Note that, for $k\ge 0$, the denominator of $M_{k+1,k}$ is a quadratic polynomial in $k$ with distinct roots. On the other hand, for the generic Jacobi matrix $L$ of the sequences in $\cA$, we can see from~\eqref{eq:alphaev} and \eqref{eq:alphaodd} that the denominator of $\alpha_k=L_{k+1,k}$ has a double root, for every~${k\ge 0}$. Therefore, the elements of $\cC_0$ can not be obtained from elements of $\cA$ by some change of parameters.\looseness=-1

Using equation \eqref{eq:beMyz}, a simple computation gives us
\begin{gather*}
	M_{k,k}= - \dfrac{d_2}{a_2} \qquad \text{if $k$ is even}, \qquad
	M_{k,k}= \dfrac{a_2 (2 b_0+b_2) + d_2}{a_2} \qquad \text{if $k$ is odd}.
\end{gather*}

The complementary Bannai--Ito polynomials are in the subclass $\cC_0$. The matrix $M$ becomes the Jacobi matrix of the normalized Bannai--Ito polynomials using the substitution
\begin{gather*}
a_2=1,\qquad b_2=1,\qquad b_0=-1/2,\qquad b_1=-2 r_2 -1/2,\\
a_1=\rho_1+\rho_2-r_1-r_2 +1/2, \\
d_1=(1-2 r_1) \rho_1 +(2 r_2+2) \rho_2- 2 r_1+1,\qquad d_2=-2 \rho_2.
\end{gather*}
Here $r_1$, $r_2$, $\rho_1$, $\rho_2$ are the parameters used to define the coefficients of the three-term recurrence relation of the complementary Bannai--Ito polynomials in~\cite[equations (3.4) and (3.5)]{GenVZ2}. Since those recurrence coefficients have some symmetries, there are other substitutions that yield the same Jacobi matrix. For example,
\begin{gather*}
a_2=1,\qquad b_2=1,\qquad b_0=-1/2,\qquad b_1= 2 \rho_1 +1/2,\qquad a_1=\rho_1+\rho_2-r_1-r_2 +1/2, \\
 d_1=(2 r_1) r_2 +(1- 2 \rho_1) \rho_2 - r_1 +1/2,\qquad d_2=-2 \rho_2.
\end{gather*}
We have in this case $w=b_0+b_1=2 \rho_1$. This value for $w$ corresponds to the one chosen in~\cite[equation~(5.11)]{TsujiDunkl}, where the authors take $w=\rho_1$ to define the complementary Bannai--Ito polynomials. The multiplier $2$ in our $w$ corresponds to a re-scaling of the variable of the polynomials.\looseness=-1

We present next some examples of the recurrence coefficients of polynomial sequences in $\cC_0$.
The substitution
\begin{gather*}
 a_1=(2r-1) a_2/2, \qquad b_2=1,\\
 d_1=(r^2-r -(2r-1)(b_0+b_1) +1) a_2/2, \qquad d_2=-(b_0+b_1+1)a_2
 \end{gather*}
yields the Jacobi matrix $M$ given by
\begin{gather*}
	M_{2k+1,2k}= -(k+r) (k+1+b_1), \qquad
	M_{2k,2k-1}= -k (k-b_1+r-1),\\
	M_{2k,2k}= b_0+b_1+1,\qquad
	M_{2k+1,2k+1}= b_0-b_1.
\end{gather*}

With the substitution
\begin{gather*}
 a_1=a_2/2,\qquad b_2=1,\qquad d_1=0,\qquad d_2=-a_2,\qquad b_0=1,\qquad b_1=0,
\end{gather*}
we obtain the matrix given by
\begin{gather*}
	M_{2k+1,2k}= -\dfrac{(2k+1)^2}{4}, \qquad\!
	M_{2k,2k-1}= -\dfrac{(2k+1)^2}{4},\qquad\!
	M_{2k,2k}= 1,\qquad\!
	M_{2k+1,2k+1}= 2.
\end{gather*}

The substitution
\begin{gather*}
\begin{split}
&a_1=(r-2)a_2/2,\qquad b_2=1,\qquad b_0=-1/2,\qquad b_1=-1/2,\\
 & d_1=(r-3)a_2/2,\qquad d_2= 2 a_2
\end{split}
\end{gather*}
gives the matrix defined by
\begin{gather*}
	M_{2k+1,2k}= -\dfrac{(2k-1)^2 (2k+r-1)^2}{(4k+r-1)(4k+r+1)}, \qquad
	M_{2k,2k-1}= -\dfrac{4 k^2(2k+r)^2}{(4k+r-3)(4k+r-1)},\\
	M_{2k,2k}= -2,\qquad
	M_{2k+1,2k+1}=2.
\end{gather*}

Using the substitution
\[ a_1=1,\qquad a_2= 2,\qquad b_0=0,\qquad b_1=1,\qquad b_2=0,\qquad d_1=-1,\qquad d_2=0, \]
 the matrix $M$ becomes the Jacobi matrix of the Chebyshev polynomials of the first kind.
In this case, $M_{k,k}=0$ for $k\ge 0$, $M_{k+1,k}=1/4$ for $k \ge 1$, the nodes are $ x_k=(-1)^k$, and the eigenvalues are $h_k=-1+(-1)^k+ 2 k (-1)^k$, for $k\ge 0$.

\subsection[Another parametrization of the subclass C\_0]{Another parametrization of the subclass $\boldsymbol{\cC_0}$}
When we take the shift $w$ equal to $-d_2/a_2$, the sequences $y_k$ and $z_k$ also become simple rational functions of $k$, with different formulas for even and odd values of $k$. In order to avoid confusion with the corresponding sequences obtained with $w=b_0+b_1$, we will write $\tilde y_k$ and $\tilde z_k$ for the sequences with $w=-d_2/a_2$, and $\tilde M$ for the associated Jacobi matrix.

 Define the following rational functions:
 \begin{gather*}
	 \tilde y_e(k)= y_e(k), \qquad \tilde 	y_o(k)= - z_e(k),\qquad
	 \tilde 	z_e(k)=- y_o(k), \qquad \tilde z_o(k) = z_o(k).
\end{gather*}

Then the sequences $\tilde y_k$ an $\tilde z_k$ are given by
\begin{gather*}
\tilde y_{2k}= y_e(k), \qquad \tilde y_{2k+1}= - z_e(k),\qquad
 \tilde z_{2k}= - y_o(k), \qquad \tilde	z_{2 k+1}= z_o(k).
\end{gather*}

The explicit formulas for $\tilde M_{k+1,k}$ are obtained from
\begin{gather*}
 \tilde M_{2k+1,2k}= -z_e(k) z_o(k),\qquad  k\ge 0,\qquad
 \tilde M_{2k, 2k-1}= - y_e(k) y_o(k), \qquad  k\ge 1.
\end{gather*}

From \eqref{eq:beMyz}, a simple computation gives us in this case
\begin{gather*}
	\tilde M_{k,k}= b_0+b_1 \qquad \text{if $k$ is even}, \qquad
	\tilde M_{k,k}= b_0-b_1+b_2 \qquad \text{if $k$ is odd}.
\end{gather*}

The Jacobi matrix of the normalized complementary Bannai--Ito polynomials is obtained when we apply to $\tilde M$ the substitution{\samepage
\begin{gather*}
a_2=1,\qquad b_2=1,\qquad b_0=-1/2,\qquad b_1=2 \rho_2+1/2,\qquad a_1=\rho_1+\rho_2-r_1-r_2+1/2, \\
 d_1=2 r_1 r_2 -2 \rho_1 \rho_2 -r_1-r_2+\rho_1+1/2,\qquad d_2=-2 \rho_1.
\end{gather*}
Here we have $w=2 \rho_1$.}

There is some kind of duality between the cases with $w=b_0+b_1$ and those with $w=-d_2/a_2$ that we will try to clarify in what follows.
Let $\cP$ be the set of parameter vectors of the form $(a_1, a_2, b_0, b_1, b_2, d_1,d_2,w)$, with $a_2\ne 0$, and let $\Gamma$ be the map that sends an element $p$ of $\cP$ to the Jacobi matrix $M$, defined by \eqref{eq:Darb}--\eqref{eq:beMyz}. Then the subclass $\cC_0$ is the set of polynomial sequences whose Jacobi matrix is in the image under $\Gamma$ of the set $\cP_0$ of vectors in $\cP$ that have $w=b_0+ b_1$. Let $\cP_1$ be the set of vectors in $\cP$ that have $w=-d_2/a_2$.

Let $\Psi\colon \cP_0 \rightarrow \cP_1$ be the map that sends $p=(a_1, a_2, b_0, b_1, b_2, d_1,d_2, b_0+b_1)$ to $\tilde p=\smash{\bigl(\tilde a_1,\tilde a_2,\tilde b_0,} \allowbreak \smash{\tilde b_1,\tilde b_2,\tilde d_1,\tilde d_2,\tilde b_0+ \tilde b_1\bigr)}$, where the entries are defined by
\begin{gather*}
\tilde a_1=a_1,\qquad \tilde a_2=a_2, \qquad \tilde b_0=b_0,\qquad \tilde b_1= -b_0 -d_2/a_2,\qquad \tilde b_2=b_2, \\
\tilde d_1=d_1+(b_0+b_1) a_2/2 + d_2/2,\qquad \tilde d_2=-a_2 (b_0+b_1).
\end{gather*}
It is easy to verify that the inverse map of $\Psi$ is given by
\begin{gather*}
 a_1=\tilde a_1,\qquad a_2=\tilde a_2, \qquad b_0=\tilde b_0,\qquad b_1= -\tilde b_0 -\tilde d_2/\tilde a_2,\qquad b_2=\tilde b_2, \\
d_1=\tilde d_1+\bigl(\tilde b_0+\tilde b_1\bigr) \tilde a_2/2 + \tilde d_2/2,\qquad d_2=-\tilde a_2 (\tilde b_0+\tilde b_1),
\end{gather*}
and hence $\Psi$ is a bijective map. Therefore, every Jacobi matrix that is in the image under $\Gamma$ of~$\cP_0$ is also in the image under $\Gamma$ of $\cP_1$, and thus $\Gamma$ is not injective.

It is also easy to see that the restriction of $\Psi$ to the set $\cP_f$ of the vectors in $\cP_0$ that satisfy $(b_0+b_1) a_2+d_2=0$ is the identity map.
Therefore, for every vector of parameters in $\cP_f$ the corresponding Jacobi matrices $M$ and $\tilde M$ coincide.

Let us note that the polynomial sequences corresponding to vectors of parameters that are in $\cP_0$ but not in $\cP_f$ are associated with at least two sequences of nodes $x_k$ and two generalized difference-eigenvalue equations, with the same sequence $h_k$ of eigenvalues. This suggests that such sequences may satisfy two different discrete orthogonality relations. See \cite{qdiscr}.

The Jacobi matrix $J$ determined by a vector of parameters in $\cP_f$ is obtained, for example, by substitution of $d_2=-(b_0+b_1)a_2$ in either $M$ or \smash{$\tilde M$}. In this way, we obtain 
\begin{gather*}
	J_{2k,2k-1}= -\dfrac{8k(k a_2 b_2+ a_1 b_2-a_2 b_1) ( k(2k+1) a_2 b_2+(b_0+b_1+2 k b_2)a_1+d_1)}{(2 a_1+(4k-1) a_2) (2 a_1+(4k+1) a_2)},
\\
	J_{2k+1,2k}=-\dfrac{2 (2 a_1+(2k+1)a_2) (2 b_1 + ( 2k+1) b_2)}{2 a_1+(4k+1) a_2} \qquad \qquad \nonumber \\
\hphantom{J_{2k+1,2k}=}\hspace{0.25mm}{}
 \times 	\dfrac{( k(2k+1) a_2 b_2+(b_0+b_1+2 k b_2)a_1+d_1)}{2 a_1 +(4k+3) a_2}, \nonumber
\end{gather*}
and
$
	J_{k,k}= b_0+b_1$ if $k$ is even,
	$J_{k,k}= b_0-b_1+b_2$ if $k$ is odd.

Notice that the parameter $d_2$ does not appear in these equations.

\section{The Bannai--Ito algebra}\label{section7}
In this section, we find several concrete realizations of the Bannai--Ito algebra using infinite matrices as generators.
Let $a_2=1$, $b_2=1$, and $b_0=0$. Then we have $h_k=(a_1 +k) (-1)^k $ and ${x_k=(b_1+k) (-1)^k}$.
Let $B_1$ be the matrix representation of the generalized difference operator~$\mathcal D$ with respect to the Newton basis with nodes $x_k$. Let $B_2$ be the matrix representation with respect to the Newton basis of the operator of multiplication by the variable $t$. We have{\samepage
\begin{equation*}
 B_1= \begin{bmatrix} h_0 & 0 & 0 & 0 & \dots \\
   g_1 & h_1 & 0 & 0 & \dots \\
     0 & g_2  & h_2  & 0 & \dots \\
       0 & 0 & g_3  & h_3 & \dots \\
   \vdots  & \vdots  & \vdots  & \ddots  &  \ddots  \end{bmatrix}, \qquad
 B_2 = \begin{bmatrix} x_0 & 1 & 0 & 0 & \dots \\
       0 & x_1 & 1 & 0 & \dots \\
      0 & 0 & x_2  & 1 & \dots \\
      0 & 0 & 0 & x_3 & \ddots \\
 \vdots  & \vdots  & \vdots  & \vdots  &  \ddots  \end{bmatrix}.
\end{equation*}
	Let us recall that $g_k=x_{k-1}(h_k-h_0) + e_k$, and
	$e_k=-d_1 + (d_1+ k d_2) (-1)^k$, for $k\ge 0$.}

Define the constants
\begin{gather}
	w_1=2 d_1 +d_2 (1-2 b_1)-a_1, \nonumber\\
	 w_2=2 d_1 + d_2 (1- 2 a_1)+ a_1 ( 2 b_1 - 2 a_1), \nonumber\\
	 w_3= a_1 -b_1 - 2 d_1 +1/2,\label{eq:defwj}
\end{gather}
and let $B_3=\{B_1, B_2\} - w_3 I$, where $I$ is the infinite identity matrix and
$\{B_1,B_2\}=B_1 B_2 + B_2 B_1$ is the anti-commutator.
By straightforward computations, we can see that
\begin{gather*}
 \{B_2, B_3\}=B_1 + w_1 I, \qquad  \{B_3,B_1\}= B_2 + w_2 I, \\
 Q=B_1^2+ B_2^2+B_3^2= \bigl((a_1+d_2)^2+(a_1-b_1)(1+ a_1 - b_1) - 2 d_1 +1/4\bigr) I,
 \end{gather*}
and that $Q$ commutes with $B_1$, $B_2$ and $B_3$.

	Therefore, $B_1$, $B_2$, $B_3$ are generators of the Bannai--Ito algebra, with structure constants $w_1$,~$w_2$,~$w_3$.

We obtain another realization of the Bannai--Ito algebra as follows. Let
\begin{equation*}
 L_1= \begin{bmatrix} h_0 & 0 & 0 & 0 & \dots \\
   0  & h_1 & 0 & 0 & \dots \\
     0 & 0  & h_2  & 0 & \dots \\
       0 & 0 & 0  & h_3 & \dots \\
   \vdots  & \vdots  & \vdots  & \vdots  &  \ddots  \end{bmatrix}, \qquad
 L_2 = \begin{bmatrix} \beta_0 & 1 & 0 & 0 & \dots \\
       \alpha_1 & \beta_1 & 1 & 0 & \dots \\
      0 & \alpha_2 & \beta_2  & 1 & \dots \\
      0 & 0 & \alpha_3 & \beta_3 & \dots \\
 \vdots  & \vdots  & \vdots  & \vdots  &  \ddots  \end{bmatrix}, \end{equation*}
	where $\alpha_n$ and $\beta_n$ are defined in \eqref{eq:alphaev}--\eqref{eq:betaodd},
	and define $L_3=\{L_1,L_2\} - w_3 I$.

	Since we have $L_1=C B_1 C^{-1}$ and $L_2= C B_2 C^{-1}$, where $C$ is the matrix defined in \eqref{eq:cnk},
 it is easy to verify that $L_1$, $L_2$, $L_3$ are generators of the Bannai--Ito algebra, with the same structure constants $w_1$, $w_2$, $w_3$ of the generators $B_1$, $B_2$, $B_3$. Let us note that $L_2$ is a Jacobi matrix.

	We can give particular values to the parameters in \eqref{eq:defwj} to obtain simpler matrices and structure constants, for example, if we substitute
	\[ a_1=1/2,\qquad a_2=1,\qquad b_1=1/2,\qquad b_2=1,\qquad d_1=1/4,\qquad d_2=-1/2 \]
in $L_1$, $L_2$, $L_3$, we obtain the matrices
	\begin{gather*}
		\tilde{L}_1= \begin{bmatrix} 1/2 & 0 & 0 & 0 & 0 & \dots \\
			0  & -3/2 & 0 & 0 & 0 & \dots \\
			0 & 0  & 5/2  & 0 & 0  & \dots \\
			0 & 0 & 0  & -7/2 & 0 & \dots \\
      0 & 0 & 0  & 0  & 9/2 & \dots \\
		\vdots  & \vdots  & \vdots  & \vdots & \vdots  &  \ddots  \end{bmatrix},
\\
	\tilde{L}_2 = \begin{bmatrix} 0 & 1 & 0 & 0 & 0 & \dots \\
			-1/4 & 0 & 1 & 0 & 0 & \dots \\
			0 & -1 & 0  & 1 &0 & \dots \\
			0 & 0 & -9/4 & 0 & 1 & \dots \\
      0 & 0 & 0  & -4  & 0  & \dots \\
		\vdots  & \vdots  & \vdots  & \vdots  & \vdots  & \ddots  \end{bmatrix},
\qquad
		\tilde{L}_3 = \begin{bmatrix} 0 & -1 & 0 & 0 & 0 & \dots \\
			1/4 & 0 & 1 & 0 & 0 & \dots \\
			0 & -1 & 0  & -1 &0 & \dots \\
			0 & 0 & 9/4 & 0 & 1 & \dots \\
      0 & 0 & 0  & -4  & 0  & \dots \\
		\vdots  & \vdots  & \vdots  & \vdots  & \vdots  & \ddots  \end{bmatrix}.
\end{gather*}
These matrices satisfy{\samepage
\[ \bigl\{\tilde{L}_1, \tilde{L}_2\bigr\}= \tilde{L}_3,\qquad \bigl\{\tilde{L}_2, \tilde{L}_3\bigr\}= \tilde{L}_1, \qquad \bigl\{\tilde{L}_3, \tilde{L}_1\bigr\}= \tilde{L}_2, \]
and the structure constants are equal to zero.}

By substitution of
		\[ a_1=1/2,\qquad a_2=1,\qquad b_1=1/2,\qquad b_2=1,\qquad d_1=1/4,\qquad d_2=-1/2 \]
in $B_1$, $B_2$, $B_3$, we obtain the matrices
	\begin{gather*}
		\tilde{B}_1= \begin{bmatrix} 1/2 & 0 & 0 & 0 & 0 & \dots \\
			-1  & -3/2 & 0 & 0 & 0 & \dots \\
			0 & -4  & 5/2  & 0 & 0  & \dots \\
			0 & 0 & -9  & -7/2 & 0 & \dots \\
      0 & 0 & 0  & -16  & 9/2 & \dots \\
		\vdots  & \vdots  & \vdots  & \vdots & \vdots  &  \ddots  \end{bmatrix},
\\
	\tilde{B}_2 = \begin{bmatrix} 1/2 & 1 & 0 & 0 & 0 & \dots \\
			0 & -3/2 & 1 & 0 & 0 & \dots \\
			0 & 0  & 5/2  & 1 &0 & \dots \\
			0 & 0 & 0 & -7/2 & 1 & \dots \\
      0 & 0 & 0  & 0  & 9/2  & \dots \\
		\vdots  & \vdots  & \vdots  & \vdots  & \vdots  & \ddots  \end{bmatrix},
\\
		\tilde{B}_3 = \begin{bmatrix} -1/2 & -1 & 0 & 0 & 0 & \dots \\
			1 & -1/2 & 1 & 0 & 0 & \dots \\
			0 & -4 & -1/2  & -1 &0 & \dots \\
			0 & 0 & 9 & -1/2 & 1 & \dots \\
      0 & 0 & 0  & -16  & -1/2  & \dots \\
		\vdots  & \vdots  & \vdots  & \vdots  & \vdots  & \ddots  \end{bmatrix}.
\end{gather*}
These matrices satisfy
\[ \bigl\{\tilde{B}_1, \tilde{B}_2\bigr\}= \tilde{B}_3,\qquad \bigl\{\tilde{B}_2, \tilde{B}_3\bigr\}= \tilde{B}_1, \qquad \bigl\{\tilde{B}_3, \tilde{B}_1\bigr\}= \tilde{B}_2, \]
			and hence the structure constants are equal to zero. We also have \smash{$\tilde{B}_1^2+ \tilde{B}_2^2+\tilde{B}_3^2= (1/4) I$}.

The Bannai--Ito algebra with zero structure constants has been extensively studied because its commutation relations can be considered as an anticommutator analog of the ordinary spin algebra $\mathfrak{su}(2)$.
See \cite{Brown} and \cite{GVYZ} for details about the finite-dimensional representations and the applications of this algebra.
			
\section{Final remarks}\label{section8}
The polynomial sequences in the class of complementary $-1$ polynomials that are not in the subclass $\cC_0$ have complicated recurrence coefficients and their study requires further research work.

The discrete orthogonality of the $-1$ polynomials associated with a sequence of pairwise distinct nodes can be studied using the same approach that we used in \cite{qdiscr}.

The construction for the $-1$ polynomials of a scheme analogous to the ones obtained by Koornwinder in \cite{KAsk,K2,K3} for the $q=1$ and the general $q$-polynomials is an interesting research project that requires further work.

Another interesting problem is that of finding a systematic way, without using limits, to associate to a given family of $q$-hypergeometric polynomials a family of $-1$ polynomials that coincides with one that can be obtained by taking limits as $q$ goes to $-1$ of the given family of $q$-polynomials.

 Finding the polynomial sequences in the class $\cC$ for which the generalized difference operator~$\cD$ can be interpreted as some kind of simple modified difference or differential operator, or some sort of Dunkl type operator, such as the operator obtained in \cite{GenVZ2} for the complementary Bannai--Ito polynomials, is also an interesting problem.

The study of how our approach compares with the one of Terwilliger may lead to some generalizations and simplifications of the results already obtained.

\subsection*{Acknowledgements}
I am very grateful to the anonymous referees that contributed numerous ideas that were used to improve the paper.

\pdfbookmark[1]{References}{ref}
\LastPageEnding

\end{document}